\documentclass{article}

\usepackage[margin=1in]{geometry} 
\usepackage{amsmath,amsthm,amssymb}
\usepackage{graphicx}
\usepackage{color}
\usepackage[dvipsnames]{xcolor}
\usepackage{subcaption}
\usepackage[export]{adjustbox}
\usepackage{url}

\newcommand{\norm}[1]{\left\|#1\right\|}

\graphicspath{ {./nlmc/} }

\DeclareMathOperator*{\argmin}{\text{argmin}}

\begin{document}

\title{Deep Multiscale Model Learning}
\author{
 Yating Wang\thanks{Department of Mathematics, Texas A\&M University, College Station, TX 77843, USA (\texttt{wytgloria@math.tamu.edu})}
\and
Siu Wun Cheung\thanks{Department of Mathematics, Texas A\&M University, College Station, TX 77843, USA (\texttt{tonycsw2905@math.tamu.edu})}
\and
Eric T. Chung\thanks{Department of Mathematics, The Chinese University of Hong Kong, Shatin, New Territories, Hong Kong SAR, China (\texttt{tschung@math.cuhk.edu.hk}) }
\and
Yalchin Efendiev\thanks{Department of Mathematics \& Institute for Scientific Computation (ISC), Texas A\&M University,
College Station, Texas, USA (\texttt{efendiev@math.tamu.edu})}
\and
Min Wang\thanks{Department of Mathematics, Texas A\&M University, College Station, TX 77843, USA (\texttt{wangmin@math.tamu.edu})}
}

\maketitle

\begin{abstract}

The objective of this paper is to design novel multi-layer neural 
network architectures for multiscale simulations of  flows 
taking into account 
the observed data and physical modeling concepts. Our approaches 
use deep learning concepts combined with local
multiscale model reduction 
methodologies to predict flow dynamics. 
Using reduced-order model concepts is important for constructing 
robust deep learning architectures since the reduced-order models  provide 
fewer degrees of freedom. Flow dynamics can be thought of as multi-layer 
networks. More precisely, the solution (e.g., pressures and saturations)
at the time instant $n+1$ depends on the solution at the time instant $n$ 
and input parameters, such as permeability fields, forcing terms,
and initial conditions. One can regard the solution 
as a multi-layer network, where each layer, in general,
 is a nonlinear forward map and the number of layers relates to the internal
time steps.
We will rely on rigorous model reduction concepts to define unknowns and 
connections for each layer.
 In each layer, our reduced-order models will provide a forward map, 
which will be modified (``trained'') using available data. It is critical 
to use reduced-order models for this purpose, which will identify the 
regions of influence and the appropriate number of variables. 
Because of the lack of available data, the training will be 
supplemented with computational data as needed and the interpolation 
between data-rich and data-deficient models. 
We will also use deep 
learning algorithms to train the elements of the 
reduced model discrete system. 
We will present main ingredients of our approach and numerical results.
Numerical results show that using deep learning and multiscale models,
we can improve the forward models, which are conditioned to the available
data.

\end{abstract}

\section{Introduction}



Many processes have multiple scales and uncertainties at the finest
scales. These include, for example, porous media processes, 
where the media properties can vary over many scales.
Constructing models on a computational coarse grid
is challenging. Many multiscale methods \cite{ egw10, GMsFEM13, lj06_2, hfmq98,ee03} and solvers are designed
to construct 
coarse spaces and resolve unresolved scales to a desired accuracy 
via additional computing. 
In general, for nonlinear problems and in the presence of observed 
solution-related data,
multiscale models are challenging to construct \cite{ehg04,cegg13,nonlinear_AM2015}. Multiscale methods,
in many challenging problems, can give a guidance to construct
robust computational models by combining multiscale concepts with
deep learning methodologies. This is an objective of this paper.


In this paper, we consider multiscale methods for nonlinear PDEs
and incorporate the data to modify the resulting coarse-grid
model. This is a typical situation in many applications, where
multiscale methods are often used to guide coarse-grid models.
These approximations, e.g., typically involve a form of the coarse-grid
equation \cite{ab05, egw10,  arbogast02, GMsFEM13, AdaptiveGMsFEM, brown2014multiscale, ElasticGMsFEM, ee03, abdul_yun, ohl12, fish2004space, fish2008mathematical, oz07, matache2002two, henning2009heterogeneous, OnlineStokes, chung2017DGstokes,WaveGMsFEM, MsDG}, 
where the coarse-grid equations are formed and the
parameters are computed or found via inverse problems \cite{inverse_Zabaras, made08, IB2013, memd10, deepconv_Zabaras}. 
 As was
shown \cite{NLMC}  the form of upscaled and multiscale equations
can be complicated, even for linear problems. To condition these
models to the available observed data, we propose a multi-layer neural network,
which uses multiscale concepts. We also discuss using deep learning techniques in approximating the coarse-grid parameters.

In this work, we will use the non-local multi-continuum approach (NLMC), developed
in \cite{NLMC}. This approach identifies the coarse-grid parameters
in each cell and their connectivity to neighboring variables. 
The approach derives its foundation
from the Constraint Energy Minimizing Generalized Multiscale Finite Element
Method (CEM-GMsFEM) \cite{chung2017constraint}, 
which has a convergence rate $H/\Lambda$, where
$\Lambda$ represents the local heterogeneities. Using the concept of
CEM-GMsFEM, NLMC defines new basis functions such that the degrees
of freedom have physical meanings (in this case, they represent
the solution averages). In this work, NLMC will be used as our 
multiscale method.


Deep learning has attracted a lot of attention in a wide class 
of applications and 
gains great success in many computer vision tasks including 
image recognition, language translation and so on \cite{dcnn_im2012, dnn_speech_2012, HK_resnet2016}. Deep Neural Network is one particular branch of artificial neural network 
algorithm under the concept of machine learning.
They are information processing systems inspired by the biological 
nervous systems and animal brains. In an artificial neural network, 
there are a collection of connected units called artificial neurons, 
which are analogous to axons in the brain of an animal or human. 
Each neuron can transmit a signal to another neuron through 
the connections. The receiving neuron will then process the 
signal and transmit the signal to downstream neurons, etc. 
Many researches have focused on learning the expressivity 
of deep neural nets theoretically
 \cite{Cybenko1989, Hornik1991, Csáji2001, Telgrasky2016, Poggio2016, Hanin2017}. 

There are numerous results to investigate the universal approximation 
property of neural networks and show the ability of deep 
networks in approximations of a rich classes of functions.
The structure of a deep neural network is usually a composition 
of multiple layers, with several neurons in each layer. 
In deep learning, each level transforms its input data into a 
little more abstract representation. In between layers, some activation 
functions are needed as the nonlinear transformation on the input 
signal to determine whether a neuron is activated or not. 
The composition structure of the deep nets is important for approximating complicated 
functions. This encourages many works utilizing deep learning 
in solving partial differential equations and model reductions. 
For example, in the work \cite{E_deepRitz} the authors numerically 
solve Poisson problems and eigenvalue problems in the context of the Ritz method based on representing the trail functions by deep 
neural networks. In \cite{Ying_paraPDE}, a neural network was proposed 
to learn the physical quantity of interest as a function of random input coefficients; the accuracy and efficiency of 
the approach for solving parametric PDE problems was shown. 
In the work \cite{deepconv_Zabaras}, the authors study
deep convolution networks for surrogate models.
In \cite{Shi_resnet}, the authors build a connection between residual 
networks (ResNet) and the characteristic equation transport equation. 
This work proposes a continuous flow model for ResNet and 
shows an alternative perspective to understand deep neural networks.


In this work, we will bring together machine learning and 
novel multiscale model reduction techniques to design/modify upscaled 
models and to train coarse-grid discrete systems. This will also allow 
alleviating some of the computational complexity involved in multiscale 
methods for time-dependent nonlinear problems. Nonlinear time-dependent 
PDEs will be treated as multi-layer networks. More precisely, the solution 
at the time instant $n+1$ depends on the solution at the time instant $n$ 
and input parameters, such as permeability fields and source terms. 
One can regard the solution as a multi-layer network. We will rely on 
rigorous multiscale concepts, for example from \cite{NLMC}, to define 
unknowns and regions of influence (oversampling neighborhood structure). 
In each layer, our reduced-order models will provide a forward map, which 
will be modified (``trained'') using available data. It is critical 
to use reduced-order models for this purpose, which will identify 
the regions of influence and the appropriate number of variables.

Because of the lack of available data in porous media applications,
the training will be supplemented with computational data as needed, 
which will result in data based modified multiscale models.  
In this work, we will consider various sources for ``real'' data, 
for example, the real data can be selected from different permeability 
fields (or can be taken as different multi-phase models), to
test our approaches. We will investigate the interpolation between the 
data-rich and data-deficient models. We will use the multi-scale 
hierarchical structure of porous media to construct neural networks 
that can both approximate the forward map in the governing non-linear 
equations and super resolve physical data to fine scales.

In our numerical example,
we will consider a model problem, a diffusion equation, and measure the 
solution at different time steps. The neural network is constructed 
using an upscaled model based on the non-local multi-continuum approach 
\cite{NLMC}. 
We have tested various neural network architectures and initializations. 
The neural network is constructed based on multiple layers.
We have selected the number of 
coarse-grid variables a priori in our simulations (based on the possible 
number of channels) and impose a constraint on the connection between 
different layers of neurons to indicate the region of influence. 
Because of the coarseness of the model, the prediction is more robust 
and computationally inexpensive.  We have observed that the network 
identifies the multiscale features of the solution and the update 
of the weight matrix correlates to the multiscale features. 

In our simulations, we train the solution using the observed data and 
computational model. The observed data is obtained from a modified
 ``true'' model with different channel permeability structure. 
We plot errors across different samples and observe that if only 
the computational model is used in the training, the error can be larger 
compared to if we use observed data in addition for the training when the 
results are close to the true model. The resulting deep neural network 
provides a modified forward map, which provides a new coarse-grid model
that is more ``accurate.'' Our approach 
indicates that incorporating some observation data in the 
training can improve the coarse grid model. The resulting deep neural 
network provides a modified forward map, which provides a new 
coarse-grid model. 
We have also observed that incorporating computational data to the existing
observed data in the training can improve the predictions, when there
is not sufficient observed data.
We have also tested deep learning algorithms for 
training elements of the stiffness matrix and multiscale basis 
functions for channelized systems. Our initial numerical results 
show that one can achieve a high accuracy using multi-layer networks 
in predicting the discrete coarse-grid systems.

The paper is organized as follows. In the next section, Section
\ref{sec:prelim}, we present general multiscale concepts.
Section \ref{sec:model} is dedicated to neural network construction.
In Section \ref{sec:numres}, we present numerical results.

\section{Preliminaries}
\label{sec:prelim}

In general, we study 
\begin{equation}
\label{eq:non1}
u_t = F(x, t, u, \nabla u, I)
\end{equation}
where $I$ denotes the input, which can
include the media properties, such as permeability field, 
source terms (well rates), or initial conditions.
$F$ can have a multiscale dependence with respect to 
space and time. The coarse-grid equation for
(\ref{eq:non1}) can have a complicated form 
for many problems (cf. \cite{NLMC}). This involves
multiple coarse-grid variables in each computational coarse
grid, non-local connectivities between the coarse-grid
variables, and complex nonlinear local
problems with constraints. In a formal way, the coarse-grid
equations in the time interval $[t_n,t_{n+1}]$ can be written
for $u_{i}^{j,n}$, where $i$ is the coarse-grid block, $j$ is
a continuum representing the coarse-grid variables, 
and $n$ is the time step. More precisely,
for each coarse-grid block $i$, one may need several
coarse-grid variables, which will be denoted by $j$.
The equation for $u_{i}^{j,n}$, in general, has a form
\begin{equation}
\label{eq:non1-discrete}
\overline{u}_i^{j,n+1} - \overline{u}_i^{j,n}  = 
\sum_{i,j}\overline{F}_{i,j}(x, t, \overline{u}_i^{j,n}, \nabla \overline{u}_i^{j,n}, I),
\end{equation}
where the sum is taken over some neighborhood cells and 
corresponding connectivity continuum.
The computation of $\overline{F}$ can be expensive and involve local nonlinear
problems with constraints. 
In many cases, researchers use general concepts from upscaling,
for example, the number of continua, the dependence of $\overline{F}$,
non-locality, to construct multiscale models.
We propose to use
the overall concept of the complex upscaled models in conjunction
with deep learning strategies to design novel data-aware coarse-grid models.
Next, we consider a specific equation.

In the paper, we consider a special
case of (\ref{eq:non1}), 
the diffusion equation in fractured media
\begin{equation}\label{eq:diffusion}
\frac{\partial u}{\partial t}- \text{div} (\kappa(x) \lambda(t, x) \nabla u) = g(t), \quad \text{in} \quad D.
\end{equation}
subject to some boundary conditions. 
Our numerical examples 
consider the zero Neumann boundary condition $\nabla u \cdot n = 0$.
Here, $D$ is the computational domain,
 $u$ is the pressure of flow, 
$g(t)$ is a time dependent source term, and $\kappa(x)$ 
is a fixed heterogeneous fractured permeability field. 
The $\lambda(t, x)$ is some given mobility which
is time dependent and represent the nonlinearities in two-phase
flow. Our approach can be applied to nonlinear equations.
As the input parameter $I$, we will consider source terms $g(t,I)$, 
which correspond to well rates.  In general, we can also consider
permeability fields as well as initial conditions as the input
parameter.
We will modify existing upscaled models using source term configurations.


\subsection{Multiscale model: Non-local multi-continuum approach}

In this section, we describe in more details  nonlocal multi-continuum
approach following \cite{NLMC}.
In our work, we consider the diffusion problem in fractured media, 
and divide the domain $D$ into the matrix region and the fractures,
where the matrix has low conductivity and the fractures are low dimensional
objects with high conductivities. That is 
\begin{equation}
D = D_m \bigoplus_i d_i D_{f,i}
\end{equation}
where $m$ and $f$ corresponds to matrix and fracture respectively, 
and $d_i$ is the aperture of fracture $D_{f,i}$. 
Denote by $\kappa(x) = \kappa_m$ the permeability 
in the matrix, and $\kappa(x) = \kappa_i$ the permeability 
in the $i$-th fracture. The permeabilities of matrix and 
fractures can differ by orders of magnitude. 

The fine-scale solution of \eqref{eq:diffusion} 
on the fine mesh $\mathcal{T}^h$ can be obtained 
using the standard finite element scheme, with backward Euler method 
for time discretization:
\begin{equation}\label{eq:fine}
\left(\frac{u_f^{n+1} - u_f ^n }{\Delta t} ,v\right) +  (\kappa \lambda^{n+1}  \nabla u_f^{n+1}, \nabla v) = (g^{n+1}, v).
\end{equation}
Here,$(\cdot,\cdot)$ denotes the $L^2$ inner product.
In the matrix form, we have
\begin{equation}
M_f u_f^{n+1} + \Delta t A_f u_f^{n+1} = \Delta t b_f + M_f u_f^{n},
\end{equation}
where $M_f$ and $A_f$ are fine scale mass and stiffness 
matrix respectively, $b_f$ is the right hand side vector.

\begin{figure}[!h]
\centering
\includegraphics[scale=0.35]{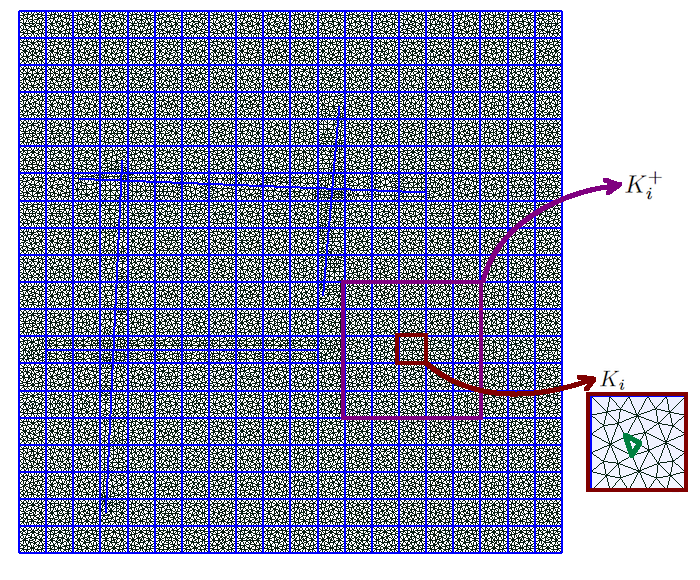}
\caption{Illustration of coarse and fine meshes. }
\label{fig:fine_coarse_mesh}
\end{figure}

For the coarse scale approximation, we assume $\mathcal{T}^H$ 
is a coarse-grid partition of the domain
$D$ with mesh size $H$
(see Figure \ref{fig:fine_coarse_mesh})
for an illustration of the fine and coarse mesh, 
where coarse elements are blue rectangles
 and fine elements are unstructured black triangles.
Denote by $\{K_i | \quad i = 1, \cdots, N\}$ the set of coarse elements
in $\mathcal{T}^H$, where $N$ is the number of coarse blocks. For each $K_i$,
we define the oversampled region $K_i^+$ to be an oversampling of $K_i$ 
with a few layers of coarse blocks. 
We will use the non-local multi-continuum approach (NLMC) \cite{NLMC}.

In the NLMC approach, 
the multiscale basis functions are selected such that the degrees
of freedom have physical meanings and correspond to average solutions.
This method derives its foundation
from Constraint Energy Minimizing Generalized Multiscale Finite Element
Method (CEM-GMsFEM) \cite{chung2017constraint}, and starts with the definition the auxiliary space.
The idea here is to use a constant as auxiliary basis for the matrix in each coarse block, and
 constants  for each separate fracture network within each coarse block. 
The simplified auxiliary space uses minimal degrees of freedom in each continua, thus one can obtain an
upscaled equation with a minimal size and the degrees
of freedom represent the averages of the solution over each continua.
 To construct the multiscale basis function for NLMC, we consider an 
oversampling region $K_i^+$ of coarse block $K_i$, 
the basis $\psi_m^{(i)}$ solves the following local 
constraint minimizing  problem on the fine grid
\begin{equation}\label{eq:basis}
\begin{aligned}
& a(\psi_m^{(i)}, v) + \sum_{K_j \subset K_i^+} \left(\mu_0^{(j)} \int_{K_j}v  + \sum_{m \leq L_j} \mu_m^{(j)} \int_{f_m^{(j)}} v \right) = 0, \quad \forall v\in V_0(K_i^+), \\
& \int_{K_j}\psi_m^{(i)}   = \delta_{ij} \delta_{0m}, \quad \forall K_j \subset K_i^+, \\
& \int_{f_m^{(j)}} \psi_m^{(i)}   = \delta_{ij} \delta_{nm}, \quad \forall f_m^{(j)} \in F^{(j)}, \; \forall K_j \subset K_i^+.
\end{aligned}
\end{equation}
where $ a(u, v)  = \int_{D_m } \kappa_m \lambda \nabla u \cdot \nabla v + \sum_i \int_{D_{f,i}} \kappa_i \lambda \nabla_f u \cdot \nabla_f v$. By this way of construction, the average of the basis $\phi_0^{(i)}$ equals $1$ in the matrix part of coarse element $K_i$, and equals $0$ in other coarse blocks $K_j \subset K_i^+$ as well as any fracture inside $K_i^+$. As for $ \phi_l^{(i)}$, it has average $1$ on the $l$-th fracture continua inside the coarse element $K_i$, and average $0$ in other fracture continua as well as the matrix continua of any coarse block $K_j \subset K_i^+$. It indicates that the basis functions separate the matrix and fractures, and each basis represents a continuum. 

We then define the transmissibility matrix $T$ by
\begin{equation}\label{eq:Trans}
T_{mn}^{(i,j)} = a(\psi_m^{(i)} ,\psi_n^{(j)} ).
\end{equation}
We note that $m,n$ denotes different continua, and $i,j$ are the indices 
for coarse blocks. Since the multiscale basis are constructed in oversampled 
regions, the support of multiscale basis for different coarse degrees of 
freedom will overlap, and this results in non-local transfer and effective 
properties for multi-continuum. The mass transfer between continua 
$m$ in coarse block $i$ and continua $n$ in coarse block $j$ is 
$T_{mn}^{(i,j)} ([u_T]_n^{(j)} -[u_T]_m^{(i)})$, where $[u_T]$ is the coarse 
scale solution.

With a simple index, we can write $T$ (tranmissibilities)
 in the following form
\begin{equation}
\begin{bmatrix}
    t_{11} & t_{12} & \dots  & t_{1n} \\
    t_{21} & t_{22} & \dots  & t_{2n} \\
    \vdots & \vdots & \ddots & \vdots \\
    t_{n1} & t_{n2} & \dots  & t_{nn}
\end{bmatrix}
\end{equation}
where $n = \sum_{i=1}^N (1+L_i)$, and $1+L_i$ means the one matrix continua 
plus the number of discrete fractures in coarse block $K_i$, and $N$ is 
the number of coarse blocks.

The upscaled model for the diffusion problem \eqref{eq:diffusion} will 
be as follows
\begin{equation}
M_T {u}^{n+1} + \Delta t A_T {u}^{n+1} = \Delta t b_T + M_T {u}^{n},
\end{equation}
where $A_T$ is the NLMC coarse scale transmissibility matrix, i.e.
\[
\begin{pmatrix}
   -\sum_j  t_{1j} & t_{12} & \dots  & t_{1n} \\
    t_{21} &  -\sum_j  t_{2j} & \dots  & t_{2n} \\
    \vdots & \vdots & \ddots & \vdots \\
    t_{n1} & t_{n2} & \dots  &  -\sum_j  t_{nj}
\end{pmatrix}
\]
 and $M_T$ is an approximation of coarse scale mass matrix. We note that both $A_T$ and $M_T$ are
non-local and defined for each continua. 

To this point, we obtain an upscaled model from the NLMC method. We remark that the results in \cite{NLMC} indicate that the upscaled equation in our modified method can use small local regions. 
 
\section{Deep Multiscale Model Learning (DMML)}
\label{sec:model}
\subsection{Main Idea}

We will utilize rigorous NLMC model as stated in previous section to solve the coarse scale problems and use the resulting solutions in deep learning framework to approximate $F$ in \eqref{eq:non1}. The advantages of NLMC approach lie 
in that, one can not only get very accurate approximations compared to the 
reference fine grid solutions, but the coarse grid solutions also have 
important physical meanings. That is, the coarse grid parameters are 
the average pressure in the corresponding matrix or fracture in a 
coarse block. Usually $\overline{F}$ is difficult to compute and conditioned to data. The idea of this work is to use the coarse grid information and available real data in combination with deep learning techniques to overcome this difficulty.

It's clear that the solution at the time instant $n+1$ depends on the solution at the time instant $n$ and input parameters, such as permeability/geometry of the fractured media and source terms. Here, we would like to learn the relationship of the solutions between two consecutive time instants by a multi-layer network. If we simply take only computational data in the training process, the neural network will provide a forward map to approximate our reduced-order models.

To be specific, let $m$ be the number of samples in the training set. 
Suppose for a given set of various input parameters, we use NLMC method 
to solve the problem and obtained the coarse grid solutions 
\[ \{u_{1}^{1},\cdots,  u_{1}^{n+1},\\
 u_{ 2}^{1},\cdots,  u_{2}^{n+1},\\
  \cdots, \cdots, \\
 u_{m}^{1},\cdots,  u_{m}^{n+1}\}
\]
at all time steps for these $m$ samples. Our goal is to use deep learning to train the coarse grid solutions and find a network $\mathcal{N}$ to describe the pushforward map between $u^{n}$ and $u^{n+1}$ for any training sample.
\begin{equation}\label{eq:universal}
u^{n+1} 	\sim \mathcal{N} (u^n, I^n),
\end{equation}
where $I^n$ is some input parameter which can also change with respect to time, and $\mathcal{N}$ is a multi-layer network to be trained.

\textit{Remark: } The proposed framework also includes nonlinear elliptic PDEs, where the map $\mathcal{N}$ corresponds to the linearised equation.

In deep network, we call $u^n$ and $I^n$ the input, and $u^{n+1}$ the output. One can take the coarse solutions from time step $1$ to time step $n$ as input, and from time $2$ to $n+1$ as output in the training process. In this case, a universal neural net $\mathcal{N}$ can be obtained. With that being said, the solution at time $1$ can be forwarded all the way to time $n+1$ by repeatedly applying the universal network $n$ times, that is
\begin{equation}\label{eq:multi-time1}
u^{n+1} 	\sim  \mathcal{N} (\mathcal{N} \cdots \mathcal{N} (u^1, I^1) \cdots , I^{n-1}), I^n).
\end{equation}

Then in the future testing/predicting procedure, given a new coarse scale solution at initial time $u_{\text{new}}^1$, we can also easily obtain the solution at final time step by the deep neural network 
\begin{equation}\label{eq:multi-time2}
u_{\text{new}}^{n+1} 	\sim  \mathcal{N} (\mathcal{N} \cdots \mathcal{N} (u_{\text{new}}^1, I^1) \cdots , I^{n-1}), I^n).
\end{equation}

One can also train each forward map for any two consecutive time instants as needed. That is, we will have $u^{j+1} \sim  \mathcal{N}_ j(u^j , I^j)$, for $j =1, \cdots, n $.  In this case, to predict the final time solution ${u_{\text{new}}^{n+1}}$ given the initial time solution ${u_{\text{new}}}^{1}$, we use $n$ different networks $\mathcal{N}_1, \cdots, \mathcal{N}_n$
\[
u_{\text{new}}^{n+1} 	\sim  \mathcal{N}_n (\mathcal{N}_{n-1} \cdots \mathcal{N}_1 (u_{\text{new}}^1, I^1) \cdots , I^{n-1}), I^n).
\]

We would like to remark that, besides the previous time step solutions, the other input parameters $I^n$ such as permeability or source terms can be different when entering the network at different time steps

\begin{figure}[!h]
\centering
\begin{subfigure}[t]{1.0\linewidth}
\centering
\includegraphics[scale=0.35]{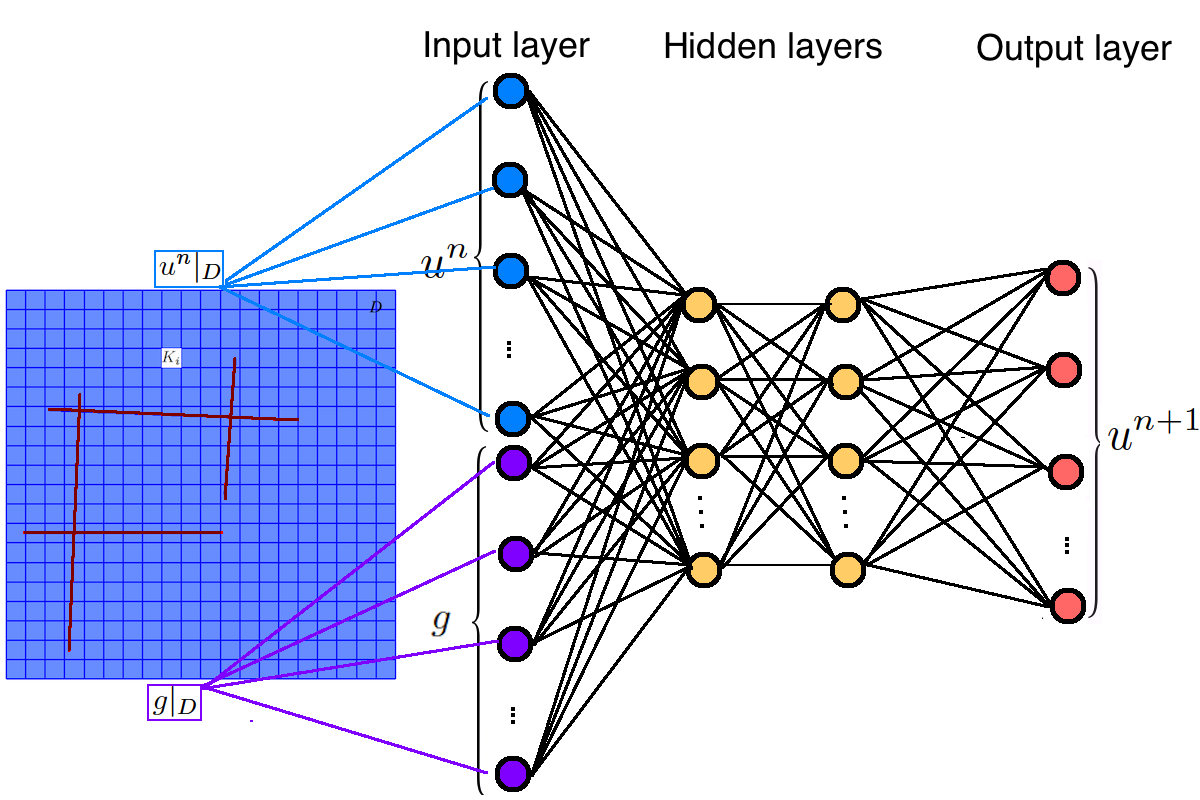}
\caption{Deep network using full input}
\end{subfigure}
\quad

\begin{subfigure}[t]{1.0\linewidth}
\centering
\includegraphics[scale=0.35]{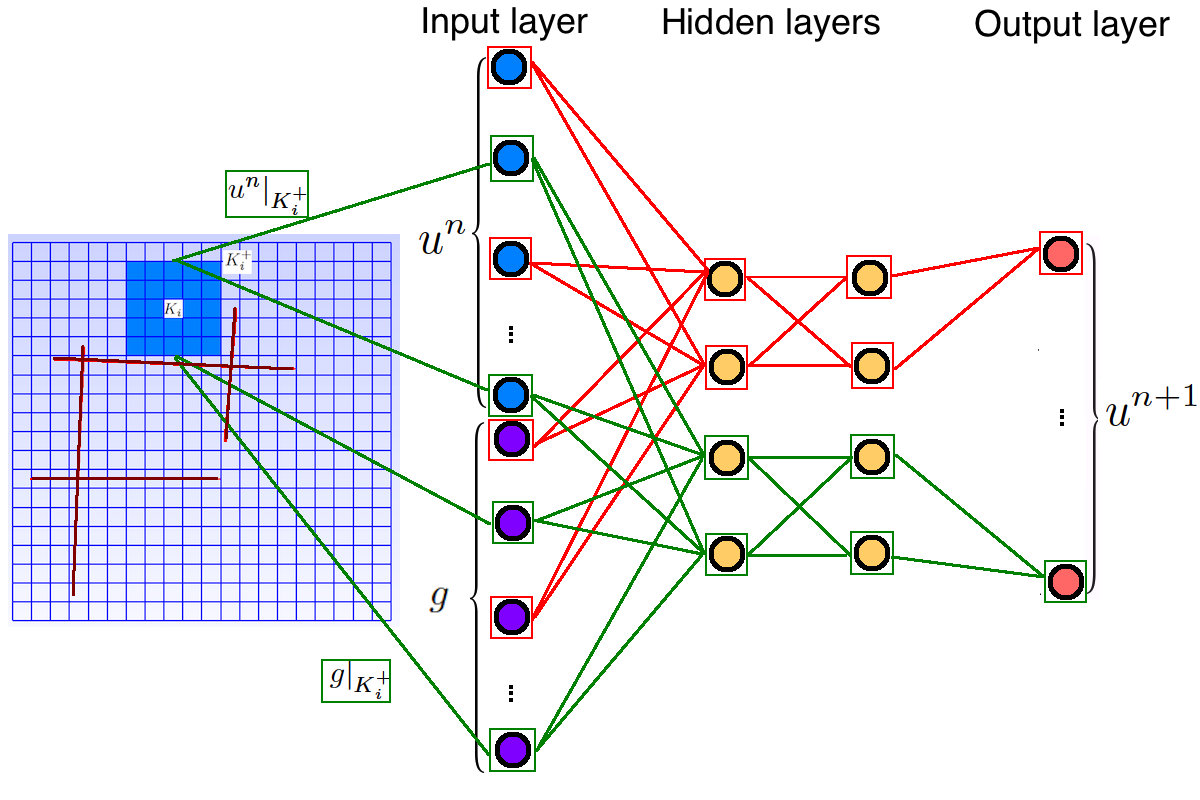}
\caption{Deep network using region of influence.}
\end{subfigure}
\caption{Comparison of deep nets with full input or with region of influence. }
\label{fig:ROI}
\end{figure}

As mentioned previously, we can also take the input in the ``region of influence''. We remark that it is important to use reduced-order model, since it will identify the regions of influence and appropriate numbers of variables. In NLMC approach, we construct a non-local multi-continuum transmissibility matrix, which provides us some information about the connections between coarse parameters. For example, for specific coarse degrees of freedom (corresponding to a coarse block or a fracture in the coarse block) of the solution at time instant $n+1$, we can simplify the problem by taking the coarse scale parameter at time instant $n$ only in the oversampling neighborhood as our input. The advantage of defining regions of the influence is to reduce the complexity of the deep network, which may also give a better initialization of the weight matrices in the training of network. An illustration of the comparison between deep neural nets with full input or with region of influence is shown in Figure \ref{fig:ROI}.

Besides all the ideas stated above, in this work, we also aim to incorporating 
available observed data in the neural net, which will modified the reduced order model and improve the performance of the model such that the new model will take into account real data effects. First, we introduce some notations.

\begin{itemize}
\item[] denote the simulation data by
\[\{ u_{s}^{1},\cdots,  u_{s}^{n+1}\}
\] 

\item[] denote the ``observation'' data by 
\[\{ u_{o}^{1},\cdots,  u_{o}^{n+1}\}
\]
\end{itemize}
at all time steps for these $m$ samples. To get the observed data, we can (1) perturb the simulation data, (2) perturb the permeability or geometry of the fractured media, run a new simulation and use the results as observed data, (3) use available experimental data. We want to investigate the effects of taking into account observation data in the output of the deep neural nets. 

 As a comparison, there are three networks we will consider:
\begin{itemize}
\item Network A: Use all observation data as output,
\begin{equation}\label{eq:No}
{u_o}^{n+1} \sim \mathcal{N}_{o} ({u_s}^n, I^n)
\end{equation}

\item Network B: Use a mixture of observation data and simulation data as output,
\begin{equation}\label{eq:Nm}
{u_{\text{mixed}}}^{n+1} \sim \mathcal{N}_m ({u_s}^n, I^n)
\end{equation}

\item Network C: Use all simulation data (no observation data) as output,
\begin{equation}\label{eq:Ns}
{u_s}^{n+1} \sim \mathcal{N}_s ({u_s}^n, I^n)
\end{equation}

\end{itemize}
where $u_{\text{mixed}}$ is a mixture of simulation data and observed data.

In Network A, we assume the observation data is sufficient, and train the observation data at time $n+1$ as a function of the observation data at time $n$. In this case, the map fits the real data in a very good manner but will ignore the simulation model if the data are obtained without using underlying simulation model in any sense. This is usually not the case in reality, since the observation data are expensive to get and deep learning requires a large amount of data to make the training effective. In Network C, we simply take all simulation data in the training process. For this network, one will get a network describes the simulation model (in our example, the NLMC model) as best as it can but ignore the observational data effects. This network can serve as an emulator (simplified
forward map, which avoids deriving/solving coarse-grid models) to do a fast simulation. We will utilize Network A and C results as references, and investigate more about Network B.  Network B is the one where we take a combination of computational data and observational data to train. It will not only take into 
account the underlying physics but also use the real data to modify the model, thus resulting in a data-driven approach.

We expect that the proposed algorithm will provide new upscaled model
that can honor the data while it follows our general multiscale concepts.


\subsection{Network structures}

Generally, in deep learning, let the function $\mathcal{N}$ be a network of $L$ layers, 
$x$ be the input and $y$ be the corresponding output. 
We write 
\begin{equation*}
\mathcal{N}(x; \theta) = \sigma(W_L \sigma (\cdots  \sigma(W_2  \sigma(W_1 x + b_1) + b_2) \cdots   ) + b_L)
\end{equation*} 
where $\theta : = (W_1, \cdots, W_L, b_1, \cdots, b_L)$, $W$'s are the weight matrices and $b$'s are the bias vectors, and $\sigma$ is the activation function. 
Suppose we are given a collection of example pairs $(x_j, y_j)$. 
The goal is then to find $\theta^*$ by solving an optimization problem
\begin{equation*}
\theta^* = \argmin_{\theta} \frac{1}{N}\sum_{j=1}^{N} ||y_j - \mathcal{N}(x_j; \theta) ||^2_2
\end{equation*}
where $N$ is the number of the samples. We note that the function $\frac{1}{N}\sum_{j=1}^{N} ||y_j - \mathcal{N}(x_j; \theta) ||^2_2$ to be optimized is called the loss function. The key points in designing the deep neural network is to choose suitable number of layers, number of neurons in each layer, the activation function, the loss function and the optimizers for the network.

In our example, without loss of generality, we suppose that there are 
uncertainties in the injection rates $g$, i.e., the value or the position 
of the sources can vary among samples.  Suppose we have a set of different realizations of the source $\{g_1, g_2, \cdots, g_m\}$, where $m$ is a sufficiently large number, we need to run simulation based on NLMC model and 
take the solutions as data for deep learning. We can perturb the geometry 
of the fractured media by translating or rotating the fractures slightly to get observation data. 

As discussed in the previous section, we consider three different networks, 
namely $\mathcal{N}_o$, $\mathcal{N}_{m}$ and $\mathcal{N}_s$. 
For each of these networks, we take the vector $x = ({u_s}^n, g^n)$ 
containing the coarse scale solution vectors and the source term in a particular time step as the input. 
As we discussed before, we can take the input coarse scale parameters in the whole domain $D$ or in the region of influence $K^+$. 
Based on the availability of the observational data in the example pairs, 
we will define an appropriate network among \eqref{eq:No}, \eqref{eq:Nm} and \eqref{eq:Ns} accordingly.
The output $y = {u_\alpha}^{n+1}$ is taken as coarse scale solution in the next time step,
where $\alpha = o,m,s$ corresponds to the network. 
Assume for extensive ensembles of source terms, there exist corresponding both computational data $u_s$ and observation data $u_o$, we will use these data to train deep neural networks $\mathcal{N}$, such that they can approximate the functions $F$ in \eqref{eq:non1} well, with respect to the loss functions. Then for some new source term $g_{m+1}$, given the coarse scale solution at time instant $n$, we expect our networks output $\mathcal{N}({u_s}^n, g^n_{{m+1}};\theta^*)$ which is close to the real data ${u_o}^{n+1}$. 

Here, we briefly summarize the architecture of the network $\mathcal{N}_\alpha$, 
where $\alpha = o,m,s$ for three networks we defined in \eqref{eq:No}, \eqref{eq:Nm} and \eqref{eq:Ns} respectively.
\begin{itemize}
\item {Input}: $x = ({u_s}^n, g^n)$ is the vector containing the coarse scale solution vectors and the source term in a particular time step. 
\item {Output}: $y = {u_\alpha}^{n+1}$ is the coarse scale solution in the next time step.
\item {Sample pairs}: $N = mn$ example pairs of $(x_j, y_j)$ are collected, where $m$ is the number of samples of flow dynamics and $n$ is the number of time steps.
\item {Standard loss function}: $\frac{1}{N}\sum_{j=1}^{N} \|y_j - \mathcal{N}_\alpha (x_j; \theta)\|_2 ^2$.
\item {Weighted loss function}: In building a network in $\mathcal{N}_m$ by using a mixture of $N_1$ pairs of observation data $\{ (x_j, y_j) \}_{j=1}^{N_1}$ and $N_2$ pairs of observation data $\{ (x_j, y_j) \}_{j=N_1+1}^{N}$, where $N_1+N_2=N$, we may consider using weighted loss function, i.e, 
$w_1 \sum_{j=1}^{N_1} \|y_j - \mathcal{N}_m (x_j; \theta)\|_2 ^2+ w_2\sum_{j=N_1+1}^{N} \|y_j - \mathcal{N}_m (x_j; \theta)\|_2 ^2$, where $w_{1} > w_{2}$ are user-defined weights.
\item {Activation function}: The popular ReLU function (the rectified linear unit activation function) is 
a common choice for activation function in training deep neural network architectures \cite{glorot11}. 
However, in optimizing a neural network with ReLU as activation function, 
weights on neurons which do not activate initially will not be adjusted, 
resulting in slow convergence. 
Alternatively, leaky ReLU can be employed to avoid such scanarios \cite{relu}. 
\item {DNN structure}: 5-10 hidden layers with 200-300 neurons in each layer.
\item {Training Optimizer}: We use AdaMax \cite{adam}, a stochastic gradient descent (SGD) type algorithm well-suited for high-dimensional parameter space, in minimizing the loss function.
\end{itemize}

\section{Numerical examples}
\label{sec:numres}

\begin{figure}[!h]
\centering
\begin{subfigure}[t]{0.45\linewidth}
\centering
\includegraphics[scale=0.3]{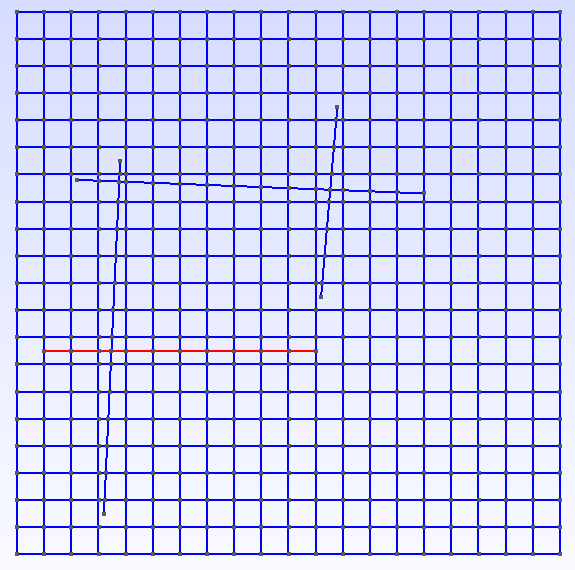}
\caption{Geometry (permeability) for obtaining simulation data.}
\label{fig:geo1}
\end{subfigure}
\begin{subfigure}[t]{0.45\linewidth}
\centering
\includegraphics[scale=0.3]{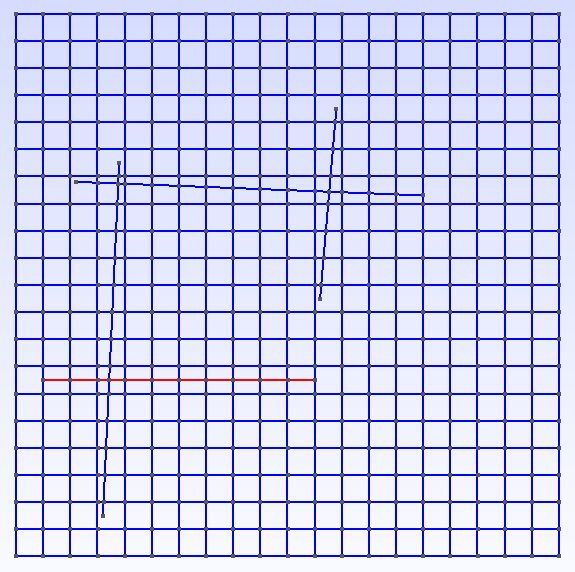}
\caption{Geometry (permeability) for obtaining observation data.}
\label{fig:geo2}
\end{subfigure}
\caption{Two geometries.}
\label{fig:2geo}
\end{figure}

In this section, we present some representative numerical results. 
In generating the NLMC model, 
we use the fractured media as shown in the Figure \ref{fig:2geo}. 
The red fracture in the two geometries are shifted up/down by one coarse 
block. To obtain the computational data, we run the simulation using the
permeability in 
Figure \ref{fig:geo1}. 
We assume that the observed data come from the solution due to the permeability
field in Figure \ref{fig:geo2}.
For the observation data,
we run
 the simulation using the permeability field on the right of Figure
 \ref{fig:geo2}. The permeability of the matrix 
is $\kappa_m = 1$, and the permeability of the fractures are $\kappa_f = 10^3$. 
We will also use a different fracture
 permeability values for the computational model.
All the network training are performed using the Python deep learning API Keras \cite{chollet2015keras}.

\subsection{Example 1}
In our first example, we use a constant mobility which is time independent. 
For the source term, we use a piecewise 
source function. Namely, in one of the coarse block, the value of $g$ is a positive number $c$, in another coarse block, the value of $g = -c$, and $g=0$ elsewhere. This is a two-well source, one of them is injection well, the other is production well, where the locations spatially change. By randomly choosing 
the location of the two wells, we get source terms $g_1, \cdots, g_{300}$. 
We run NLMC simulation for these 300 source terms on two geometries 
as shown in Figure \ref{fig:2geo}. As a result, we generate two sets of data (computational and observational data). For the $300$ source terms, we choose $290$ of them for training and $10$ for testing. We solve the equation \eqref{eq:diffusion}, set $T=0.01$, and divide it into $10$ time steps. We note that 
in this example, the value of the source is time independent. 
 
In our numerical example, we would like to find a universal deep network to describe the map between two time steps, as described in \eqref{eq:universal}. 
We use the solution at time step $1$ to time step
$9$ as input data, and from at time step $2$ to time step $10$ as output data.
Thus, the solutions corresponding to $290$ different training source terms result in $290*9 = 2610$ samples, and  the solutions corresponding to $10$ testing source terms result in $10*9 = 90$ testing samples, where the multiple $9$ is the $9$ time steps (time steps $1$ to $9$, or  time $2$ to $10$). 

We will test the performance of the three networks \eqref{eq:No}, \eqref{eq:Nm}, and \eqref{eq:Ns}. For the computational data $u_s$, we use the solution from the geometry in \ref{fig:geo1} for 300 source terms, this is the case with no real data in the training. For the observation data $u_o$, we use the solution from the geometry (permeability) in \ref{fig:geo2} for 300 source terms, this is the case with full real data in the training. As for the mixture $u_m$ of computational and observation data, we take $150$ from $u_s$ and $150$ from $u_o$, this is the case with partial real data in the training. In practice, to explain the mixture data $u_m$, we can assume we have the observation data in the whole domain given some well configurations, but for some other well configurations, we only get simulation results. In the training process, we also consider both the full input and the region of influence input (see Figure \ref{fig:ROI}), where we use multiscale
concepts to reduce the region of influence (connection) between the nodes.

The results are shown in Table \ref{tab:ex1}. First, we would like to compare the
 results between using the coarse parameters in the whole domain and 
using the coarse parameters just in the region of influence as input 
in the training. Comparing Table \ref{tab:ex1-full} and  
Table \ref{tab:ex1-sparse}, we can see that, using the region of influence idea can help to get better results for all three networks 
$\mathcal{N}_o$, $\mathcal{N}_m$ and $\mathcal{N}_s$ 
when we use similar network parameters such as the number of layers, 
number of neurons in each layer, training epochs, learning rate, loss 
functions and activation functions. This suggest that, the data in the
 region of influence can give a better initialization in the training 
compared with the data in the whole domain.

Next, we compare the results using both observation data and computational 
data, and compare the performance of the three networks defined in 
\eqref{eq:No}, \eqref{eq:Nm}, \eqref{eq:Ns}. For both sub-tables, 
we can see that, using a mixture of computational and observation 
data (the third column in the tables), we can get a better model, 
since the mean error of $\norm{\mathcal{N}_m ({u_s}^n, I^n) - {u_o}^{n+1} }$ 
among testing samples closer to the mean error of 
$\norm{\mathcal{N}_o ({u_s}^n, I^n) - {u_o}^{n+1} }$. The error history for 
some samples are also plotted. We can also observe that the deep neural 
network outputs $\mathcal{N}_m ({u_s}^n, I^n)$ (the orange curve) is closer to the observation 
data ${u_o}^{n+1}$ compared with the outputs from $\mathcal{N}_s$ (the blue curve),
 where only simulation data is used. 
We have also tested adding computational data to the observed data. In particular,
we have used only $150$ observation data and compared the results to
using  $150$ (the same) observation data and (the additional)
 $150$ computational data. The latter
provides more accurate predictions, which indicates that 
 incorporating some computational data to the observed
data can improve the predictions,  when there
is not sufficient observed data.

\begin{table}[!htb]
\centering
  \begin{tabular}{|c  |c  | c | c |}
    \hline
Errors (\%) &$\norm{\mathcal{N}_o ({u_s}^n, I^n) - {u_o}^{n+1} }$ &$\norm{\mathcal{N}_m ({u_s}^n, I^n) - {u_o}^{n+1} }$ &$\norm{\mathcal{N}_s ({u_s}^n, I^n) - {u_o}^{n+1} }$\\  \hline 
mean & 3.6  &10.6 & 19.7	\\  \hline
  \end{tabular} \caption{Example 1. Using full inputs} \label{tab:ex1-full}

  \begin{tabular}{|c  |c  | c | c |}
    \hline
Errors (\%) &$\norm{\mathcal{N}_o ({u_s}^n, I^n) - {u_o}^{n+1} }$ &$\norm{\mathcal{N}_m ({u_s}^n, I^n) - {u_o}^{n+1} }$ &$\norm{\mathcal{N}_s ({u_s}^n, I^n) - {u_o}^{n+1} }$\\  \hline 
mean &1.8 &7.5 &16.7	\\  \hline
  \end{tabular} \caption{Example 1. Using region of influence} \label{tab:ex1-sparse}
  
\caption{Example 1 (all errors are relative and in percentages). Two subtables compare the performance of deep networks when using full inputs and the region of influence in the training. In each table, three columns compare the three network using different ratio of computational data and observational data.}\label{tab:ex1}
\end{table}

\begin{figure}[!h]
\centering
\begin{subfigure}[t]{0.5\linewidth}
\centering
\includegraphics[scale=0.35]{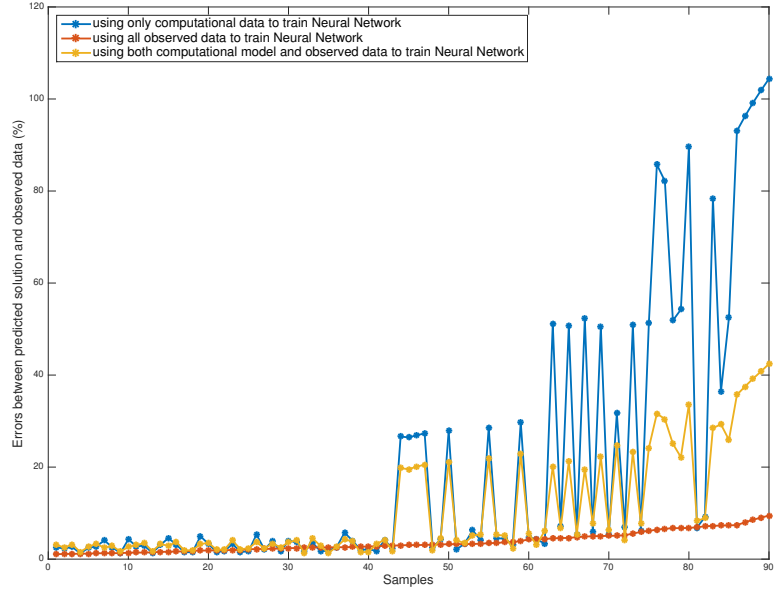}
\caption{Full input.}
\end{subfigure}

\begin{subfigure}[t]{0.5\linewidth}
\centering
\includegraphics[scale=0.35]{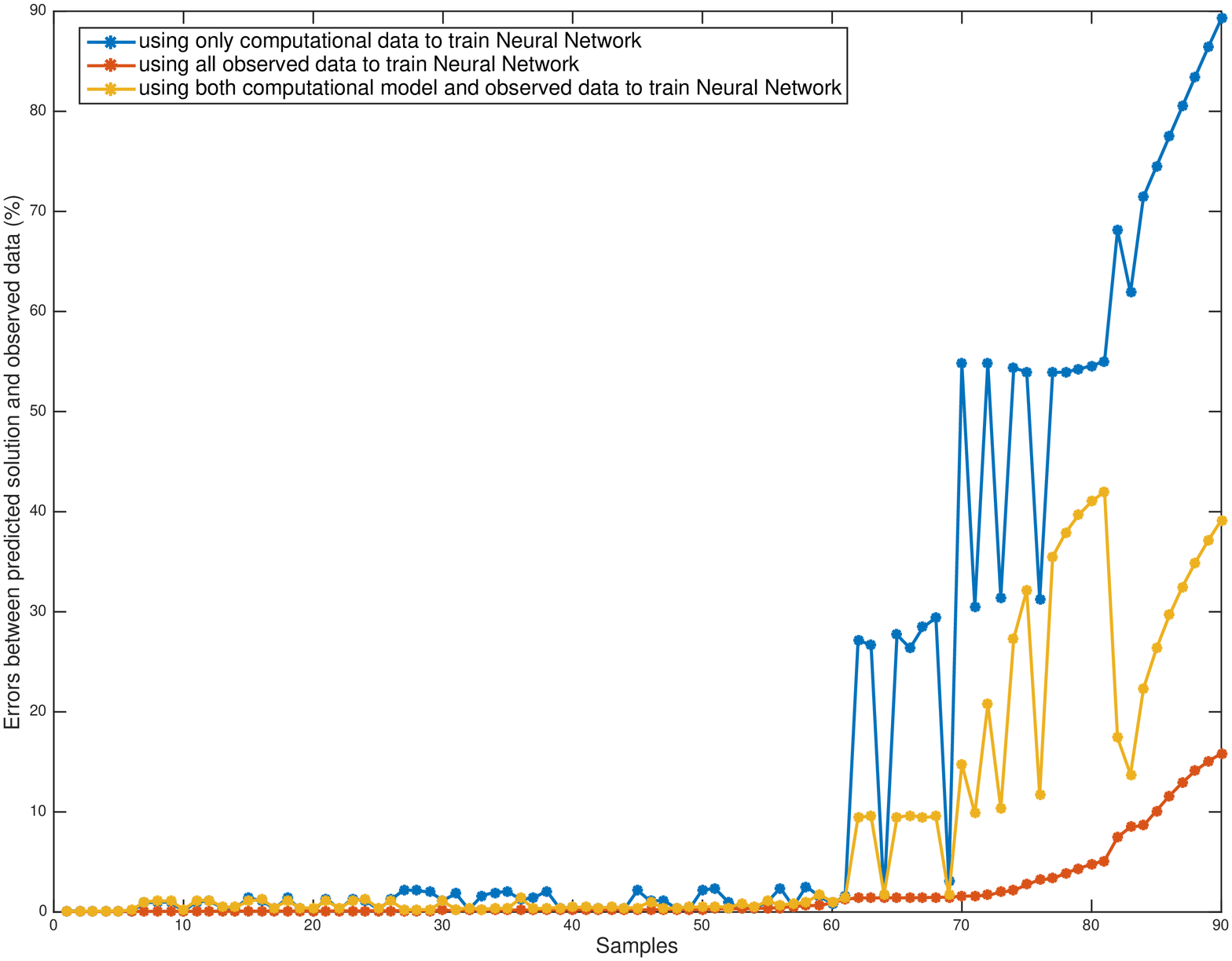}
\caption{Region of influence input.}
\end{subfigure}
\caption{Example 1. Errors of all testing samples for three networks defined in \eqref{eq:No}, \eqref{eq:Nm}, \eqref{eq:Ns}.}
\label{fig:ex1}
\end{figure}

\subsection{Example 2}

\begin{figure}[!h]
\centering
\begin{subfigure}[t]{0.3\linewidth}
\centering
\includegraphics[scale=0.35]{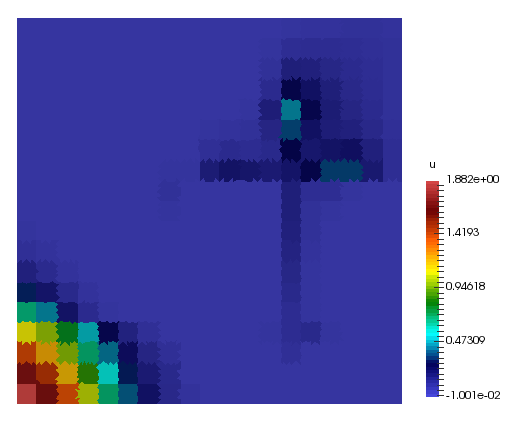}
\caption{Mobility at time $t=0.1$}
\end{subfigure}
\begin{subfigure}[t]{0.3\linewidth}
\centering
\includegraphics[scale=0.35]{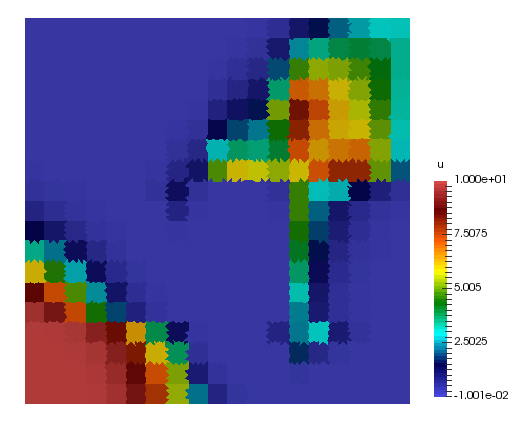}
\caption{Mobility at time $t=0.5$}
\end{subfigure}
\begin{subfigure}[t]{0.3\linewidth}
\centering
\includegraphics[scale=0.35]{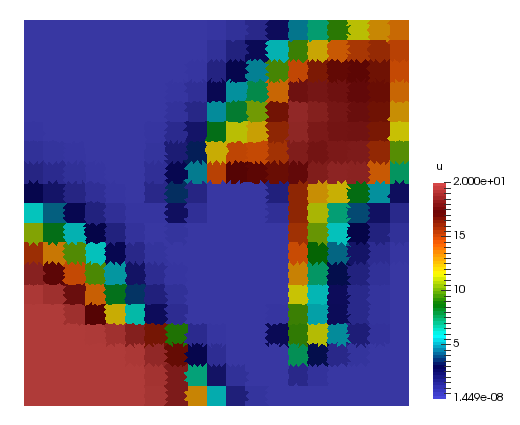}
\caption{Mobility at time $t=1.0$}
\end{subfigure}
\caption{Illustration of mobility $\lambda(x,t)$. }
\label{fig:mobility}
\end{figure}

In our second example, we use heterogeneous time-dependent mobility and 
source term. Here, we fix the location of the source term and 
vary the value of the source. The mobility is a time-dependent 
function. The distribution of the mobility in some time steps are shown in Figure \ref{fig:mobility}, which is from two-phase flow mobility.
The source term in the right hand side of the equation is piecewise 
constant functions. At $0 \leq x \leq 0.1, 0\leq y \leq 0.1$, we have 
$g=10[(\sin( \alpha x))^2+(\sin(\beta y))^2]$ denotes an injection well, 
and at $0.9 \leq x \leq 1.0, 0.9\leq y \leq 1.0$ we have 
$g=-10[(\sin(\alpha x))^2+(\sin(\beta y))^2]$ denotes an production well,
where the parameters $\alpha$ and $\beta$ are randomly chosen in each 
time step, and are different among samples (which are obtained 
using these different source terms $g$). So for each sample, we have the 
different values of the source term, and, in each sample, 
the source term is time dependent. In this example, we use $500$ different 
sources. The samples are similarly constructed as discussed in Example 1.

Again, for the computational data $u_s$, we use the solution from the 
geometry (permeability) in \ref{fig:geo1} for 500 source terms. 
For the observation data $u_o$, we use the solution from the geometry 
(permeability) in \ref{fig:geo2} for 500 source terms. In this example, 
as for the mixture $u_m$ of computational and observation data, 
we take all $500$ sample sources from $u_s$, but in half 
of the computational domain, and 
all $500$ sample sources from $u_o$ in the other half of the computational 
domain. In this example, to explain the practical meaning of $u_m$, 
we can imagine that, given all well rates, we have the observation data 
in the half domain , but in the other part of the domain we only get 
simulation results. 

\begin{figure}[!h]
\centering
\includegraphics[scale=0.35]{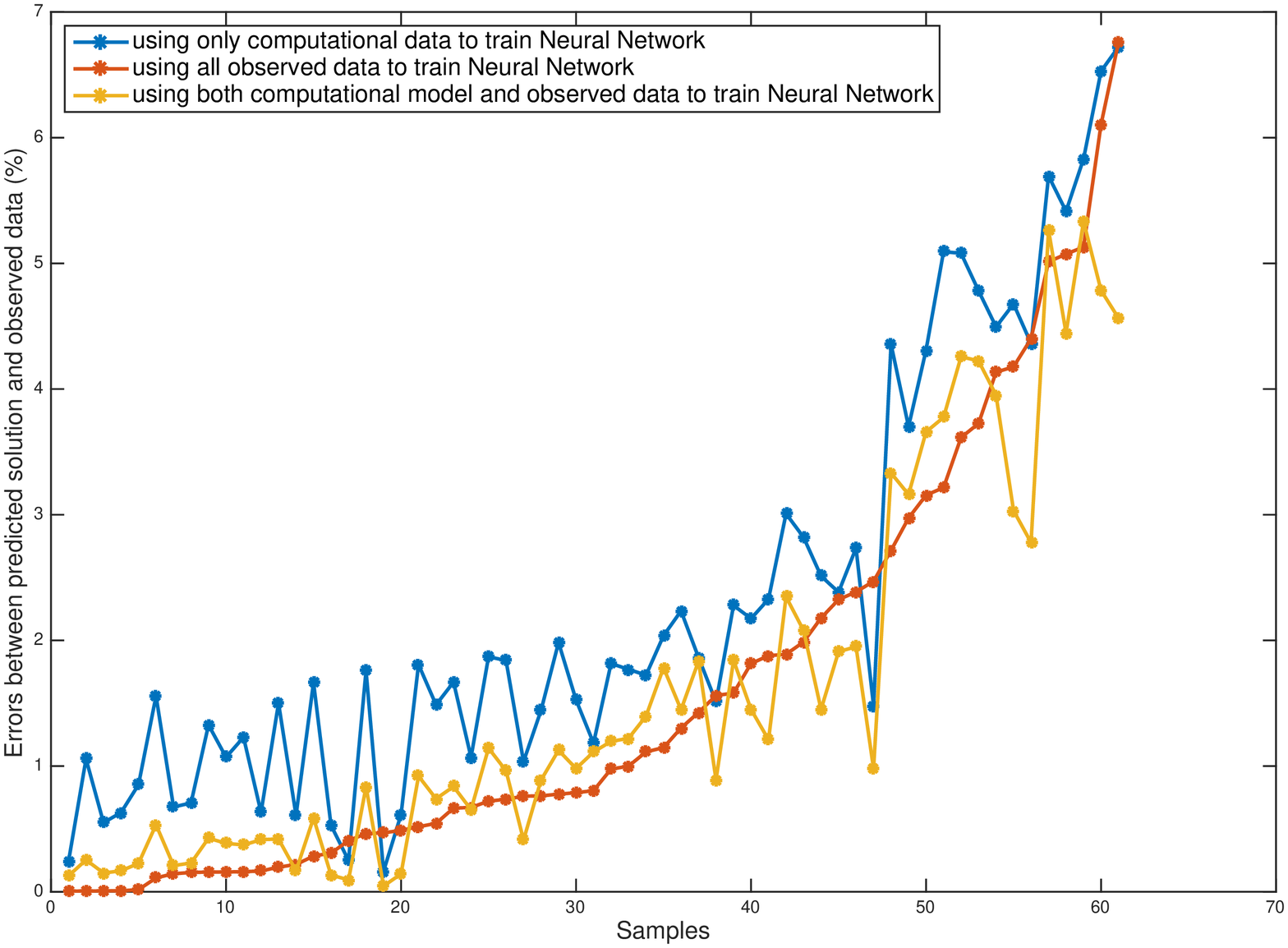}
\caption{Example 2. Using full input. Errors of all testing samples for three networks defined in \eqref{eq:No}, \eqref{eq:Nm}, \eqref{eq:Ns}.}
\label{fig:ex2}
\end{figure}

In this example, we compare the performance of the three networks 
for data-sufficient and data-deficient cases. The errors between 
the three deep networks and the real observation data are, 
$\norm{\mathcal{N}_o ({u_s}^n, I^n) - {u_o}^{n+1} }$ 
(shown in red curve in Figure \ref{fig:ex2}), 
$\norm{\mathcal{N}_m ({u_s}^n, I^n) - {u_o}^{n+1} }$
(shown in orange curve in Figure \ref{fig:ex2}), and 

$\norm{\mathcal{N}_s ({u_s}^n, I^n) - {u_o}^{n+1} }$ 
(shown in blue curve in Figure \ref{fig:ex2}) respectively. 
We use the red curve as reference, and notice that the errors 
shown in the orange curve are closer to the red one. This indicates 
that using a mixture of computational data and observation data can help 
to enhance the performance the model induced by deep learning. 
The mean errors are shown in Table \ref{err-ex2}. Although close 
(since the difference between the computation model and the 
observation model is small), it still shows the superior 
of using mixture data in the training process. In our next example,
we will change the permeability of the fracture in the computational
model, which will increase the difference between
the observed and computational data.

\begin{table}[!htb]
\centering
  \begin{tabular}{|c  |c  | c | c |}
  \hline
Errors (\%) &$\norm{\mathcal{N}_o ({u_s}^n, I^n) - {u_o}^{n+1} }$ &$\norm{\mathcal{N}_m ({u_s}^n, I^n) - {u_o}^{n+1} }$ &$\norm{\mathcal{N}_s ({u_s}^n, I^n) - {u_o}^{n+1} }$\\  \hline 
mean &1.6 &1.7 &2.2	\\  \hline
  \end{tabular}
\caption{Example 2. Mean error over testing samples using three networks.}\label{err-ex2}
\end{table} 

\begin{figure}[!h]
\centering
\includegraphics[scale=0.5]{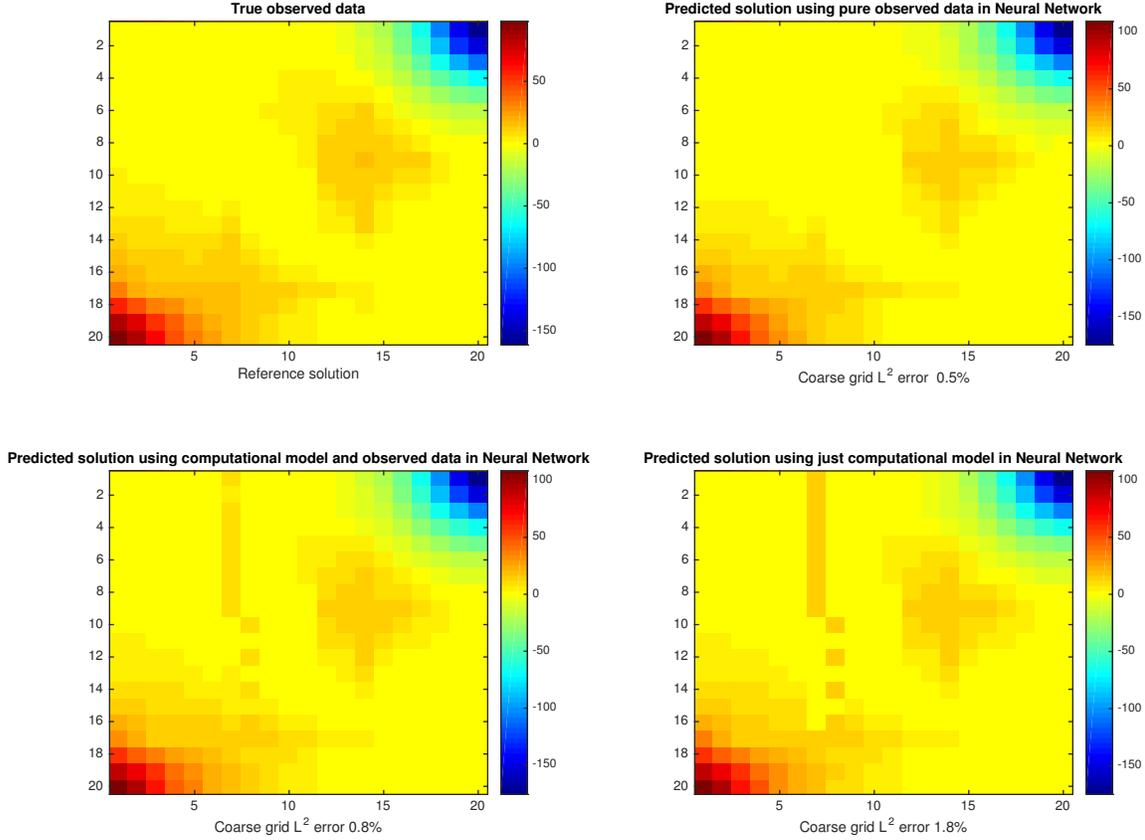}
\caption{Example 2. Predicted solutions' comparison for one of the samples.}
\label{fig:ex2-sample86}
\end{figure}

As we discussed before, we can use \eqref{eq:multi-time1} or \eqref{eq:multi-time2} to forward 
the solution from the initial time step to the final time step using 
the ``universal'' deep neural nets. Here, we will do the experiments 
using the three networks $\mathcal{N}_s$, $\mathcal{N}_m$ and 
$\mathcal{N}_s$. Actually, we assume we have 10 time steps in total, 
for given ${u_s}^1$ at the initial time, we will apply $\mathcal{N}_s$ $8$ 
times first. Then, we apply $\mathcal{N}_s$, $\mathcal{N}_m$ and 
$\mathcal{N}_s$ at the last to get the final time step predictions. That is,
\[
{u_\alpha}^{10} = \mathcal{N}_\alpha \underbrace{\mathcal{N}_s(\mathcal{N}_s \cdots ( \mathcal{N}_s}_\text{8 times}({u_s}^1))
\]
for $\alpha = o, m, s$.

\begin{figure}[!h]
\centering
\includegraphics[scale=0.3]{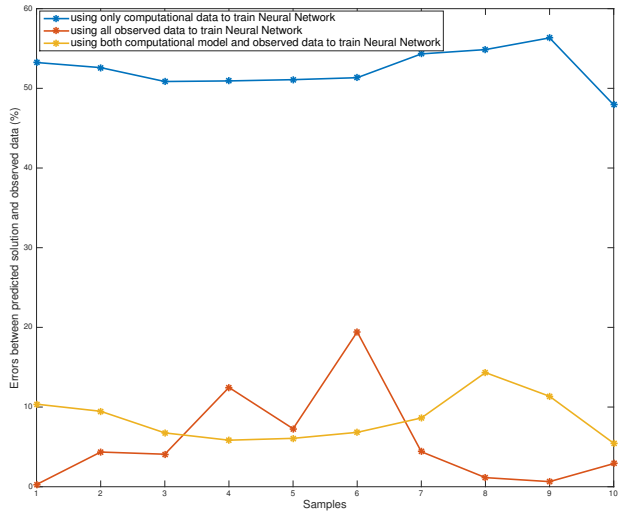}
\caption{Example 2. Errors of all 10 testing samples between predicted solutions and observed data at final time step(\%).}
\label{fig:ex2-t1-t9}
\end{figure}

Finally, we compare the final time predictions $y_\alpha$ (for $\alpha = o, m, s$) 
with the observation data at the final time step given ${u_s}^1$. 
Figure \ref{fig:ex2-t1-t9} shows the results. There are $10$ samples to test 
in total. The mean error of the red curve is $5.6\%$, the mean error of 
the orange curve is $8.5\%$, and the mean error of the blue curve is $52.3\%$. 
We can see that, the mixture-data-driven deep network predictions 
(the orange curve) and pure-observation-data-driven deep network 
predictions (the red curve) both have very good behavior as expected. 

\subsection{Example 3}

In this example, we use the geometry (permeability) shown in Figure \ref{fig:geo2}. For the observation data, we set the permeability of the fractures $k_f = 1000$ as before. For the computational data, we set the permeability of the fractures $k_f = 10$, which makes the flow within fracture is weak.
 As in Example 2, we use heterogeneous time dependent mobility and source term. Again, we fix the location of the source term and vary the value of the source. The mixture of observation data and computational data contains half of the samples from observation data, and half of the sample from computational data.
We note that for these two sets of data, the geometry stays the same, 
but the permeabilities have high contrast, thus the computational 
data are very different from the observed data. 

\begin{figure}[!h]
\centering
\includegraphics[scale=0.35]{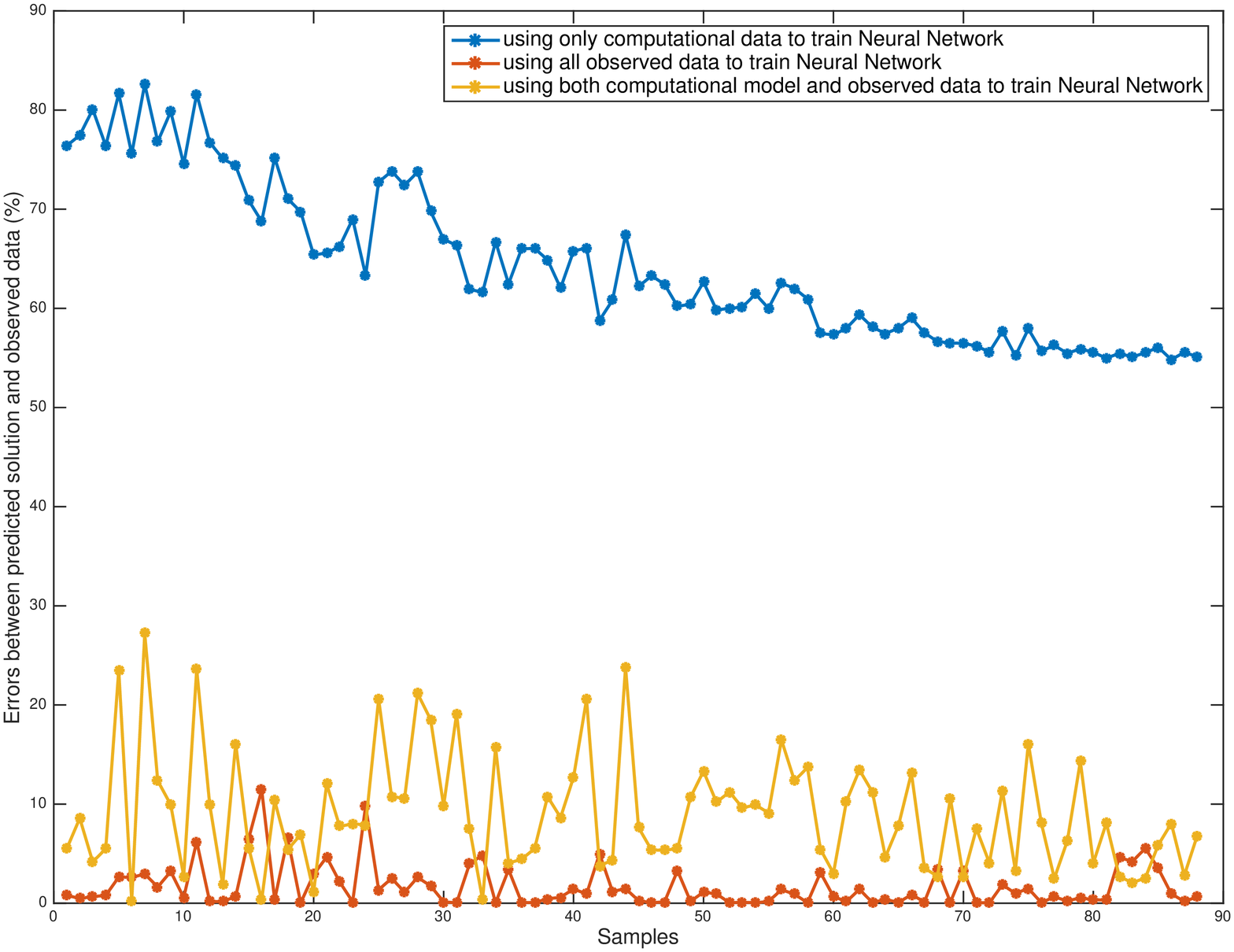}
\caption{Example 3. Using full input. Errors of all testing samples for three networks defined in \eqref{eq:No}, \eqref{eq:Nm}, \eqref{eq:Ns}.}
\label{fig:ex3}
\end{figure}

\begin{figure}[!h]
\centering
\includegraphics[scale=0.5]{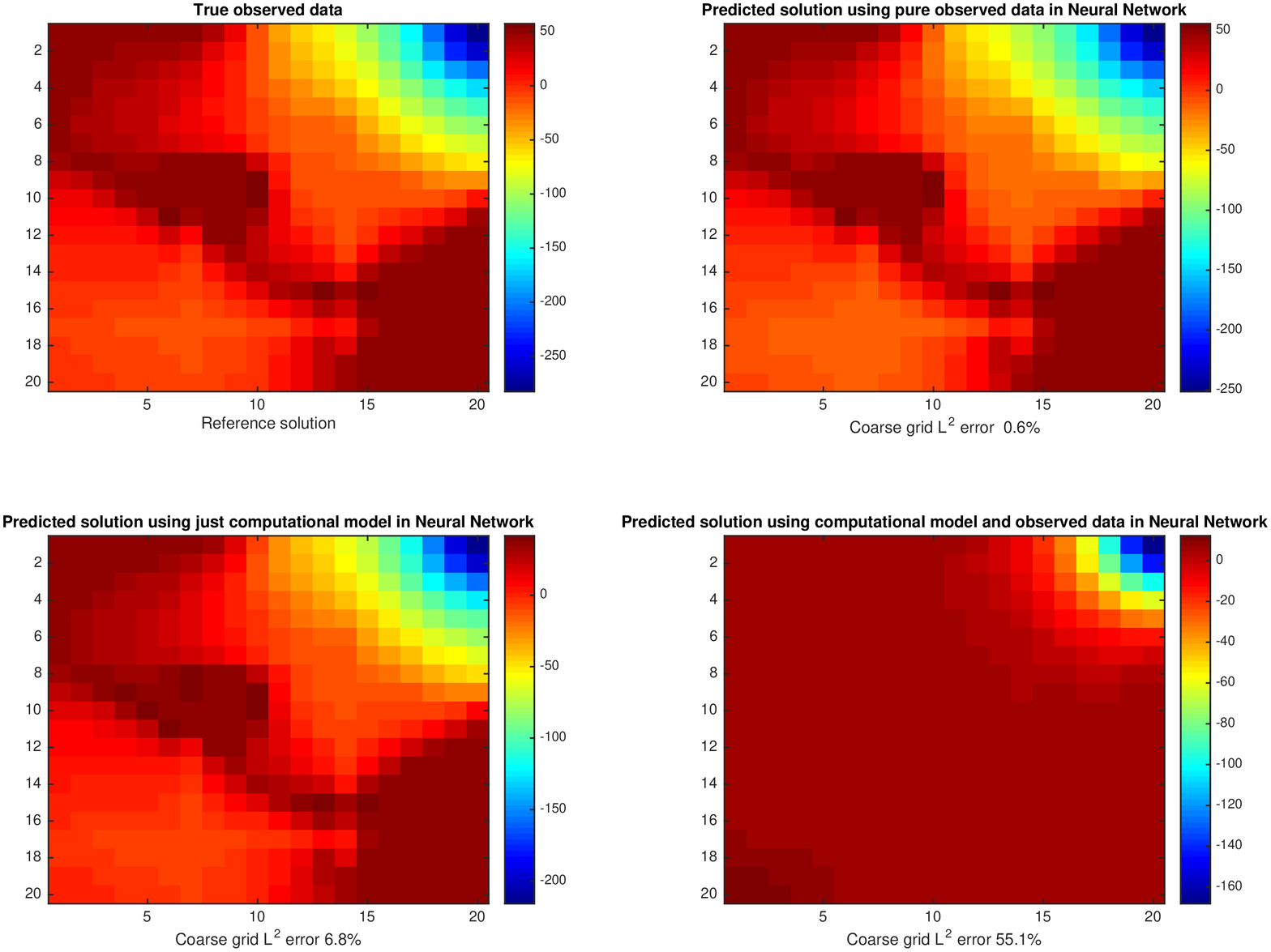}
\caption{Example 3. Predicted solutions' comparison for one of the samples.}
\label{fig:ex3-sample90}
\end{figure}

In Figure \ref{fig:ex3}, we can see that, only using the computational data in the training process is far from enough, the errors (blue curve) between the output deep network $\mathcal{N}_s ({u_s}^n, I) $  and the observation data is much larger compared with the other two curves. However, adding some observation data into the training data, the errors (orange curve) between the output deep network $\mathcal{N}_m ({u_s}^n, I) $  and the observation data is pretty good. From the Table \ref{err-ex3}, we also observe that for the mean errors across testing samples, $\norm{\mathcal{N}_m ({u_s}^n, I) - {u_o}^{n+1} }$ is much closer to $\norm{\mathcal{N}_o ({u_s}^n, I) - {u_o}^{n+1} }$. One comparison of the solutions obtained from the three networks are shown in Figure \ref{fig:ex3-sample90}, which illustrate that the network $\mathcal{N}_m ({u_s}^n, I) $ can produce reliable output.

\begin{table}[!htb]
\centering
  \begin{tabular}{|c  |c  | c | c |}
  \hline
Errors (\%) &$\norm{\mathcal{N}_o ({u_s}^n, I) - {u_o}^{n+1} }$ &$\norm{\mathcal{N}_m ({u_s}^n, I) - {u_o}^{n+1} }$ &$\norm{\mathcal{N}_s ({u_s}^n, I) - {u_o}^{n+1} }$\\  \hline 
mean &2.6 &8.8 &64.3	\\  \hline
  \end{tabular}
\caption{Example 3. Mean error over all testing samples using three networks.}\label{err-ex3}
\end{table}

\section{Conclusions}

The paper uses deep learning techniques to derive and modify
upscaled
models for nonlinear PDEs. In particular, we
combine multiscale model reduction (non-local multi-continuum
upscaling) and deep learning techniques in deriving coarse-grid
models, which take into account observed data. 
Multi-layer networks provide a nonlinear mapping between the time
steps, where the mapping has a certain structure.
The multiscale concepts, used in multi-layer networks, provide
appropriate coarse-grid variables, their connectivity information, and
some information about the mapping. However, constructing complete
and accurate nonlinear push-forward map is expensive and not possible, 
in general
multiscale simulations. Moreover, these models will not honor
the available data. In this paper, we combine the multiscale model reduction
concepts and deep learning techniques and study the use of observed
data with a new framework, Deep Multiscale Model Reduction Learning 
(DMML). We present numerical results, where we test our main concepts.
We show that the regions of influence derived from upscaling concepts
can improve the computations.
Our approach 
indicates that incorporating some observation data in the 
training can improve the coarse grid model. 
Similarly, incorporating some computational data to the observed
data can improve the predictions,  when there
is not sufficient observed data.
 The use of coarse-degrees of freedom
is another main advantage of our method. Finally, we use observed data
and show that DMML can obtain accurate solutions, which can honor the observed
data. In conclusion, we believe DMML can be used as a new coarse-grid
model for complex nonlinear problems with observed data, 
where upscaling of the computational
model is expensive and may not accurately
represent the true observed model.

\bibliographystyle{plain}
\bibliography{references}

\end{document}